\documentclass[a4paper,11pt, DIV12]{amsart}
\usepackage[english]{babel}
\usepackage[T1]{fontenc}
\usepackage[ansinew]{inputenc}
\usepackage{amsmath}
\usepackage{amsfonts}
\usepackage{ mathrsfs }
\usepackage{amssymb}
\usepackage{textcomp}
\usepackage{enumitem}
\usepackage{hyperref}
\usepackage[arrow, matrix, curve]{xy}

\theoremstyle{definition}
\newtheorem{defn}{Definition}[section]
\newtheorem{bem}[defn]{Remark}

\theoremstyle{plain}
\newtheorem{satz}[defn]{Theorem}
\newtheorem{kor}[defn]{Corollary}
\newtheorem{lem}[defn]{Lemma}

\newcommand{\A}{\mathbb{A}}
\newcommand{\Z}{\mathbb{Z}}
\newcommand{\N}{\mathbb{N}}

\newcommand{\R}{\mathbb{R}}

\newcommand{\T}{\mathbb{T}}

\newcommand{\G}{\mathbb{G}}
\newcommand{\0}{\mathcal{O}}

\newcommand{\Linear}{\mathbb{L}}

\newcommand{\CC}{\mathscr{C}}

\newcommand{\DD}{\mathscr{D}}

\newcommand{\vedge}{\land}
\newcommand{\del}{\partial}
\newcommand{\Xan}{X^{\an}}
\newcommand{\inj}{\hookrightarrow}
\newcommand{\surj}{\twoheadrightarrow}
	\DeclareMathOperator{\trop}{trop}
	\DeclareMathOperator{\Spec}{Spec}
	\DeclareMathOperator{\sing}{sing}
	\DeclareMathOperator{\an}{an}

	\DeclareMathOperator{\Trop}{Trop}
	
	\DeclareMathOperator{\pr}{pr}

	\DeclareMathOperator{\supp}{supp}
	
	\DeclareMathOperator{\id}{id}

\title[Poincar\'e lemma for differential forms on Berkovich spaces]{A Poincar\'e lemma for real-valued differential forms on Berkovich spaces}
\author{Philipp Jell}

\begin{document}
\begin{abstract}
Real-valued differential forms on Berkovich analytic spaces were introduced by Chambert-Loir and Ducros in \cite{CLD} using superforms on polyhedral complexes. We prove a Poincar\'e lemma for these superforms and use it to also prove a Poincar\'e lemma for real-valued differential forms on Berkovich spaces. For superforms we further show finite dimensionality for the associated de Rham cohomology on polyhedral complexes in all (bi-)degrees. We also show finite dimensionality for the real-valued de Rham cohomology of the analytification of an algebraic variety in some bidegrees. 
\end{abstract}

\maketitle 

\tableofcontents

\section{Introduction}

Chambert-Loir and Ducros recently introduced smooth real-valued differential forms on Berkovich analytic spaces \cite{CLD}. Basic ingredients in this new theory are Lagerberg's superforms \cite{Lagerberg} and new methods from tropical geometry (see for example \cite{Gubler2}). Chambert-Loir and Ducros also introduce currents, develop basic results of pluripotential theory (in particular an analogue of Bedford-Taylor theory), show a Poincar\'e-Lelong formula and construct Monge-Amp\`{e}re measures on Berkovich analytic spaces. There has also been recent work by Gubler and K\"unnemann \cite{GublerKuennemann} where they extend the results by Chambert-Loir and Ducros and give a Berkovich analytic construction of local heights. The results mentioned above show that the forms introduced by Chambert-Loir and Ducros are of basic interest and should play for example a crucial role in (non-Archimedean) Arakelov theory. \\
It is the aim of this paper to prove some foundational results about the differential forms of Chambert-Loir and Ducros. More precisely we will show that they satisfy a Poincar\'e lemma and we will investigate their de Rham cohomology.  \\
Our main object of study is the space of real-valued differential forms introduced by Chambert-Loir and Ducros, which are bigraded and have differential operators $d$, $d'$ and $d''$ analogous to the operators $d$, $\del$ and $\overline{\del}$ for differential forms on complex manifolds.  We are interested in the de Rham cohomology defined by these operators. Since these forms are locally defined by superforms on polyhedral complexes, we will prove a Poincar\'e lemma for superforms on polyhedral complexes. From this we can deduce finiteness of the de Rham cohomology of superforms on polyhedral complexes, using techniques analogous to those from differential geometry. We will also deduce a Poincar\'e lemma for forms on Berkovich spaces. With the help of sheaf theory we can then deduce that in some bidegrees this cohomology depends only on the underlying topological space of the Berkovich space and agrees with singular cohomology. Finally the theory of skeletons enables us to show that in these degrees cohomology is finite dimensional in many cases. \\
In section \ref{secPLemma} we recall the construction of superforms on polyhedral complexes (as introduced by Lagerberg in \cite{Lagerberg}, see as well \cite{CLD} and \cite{Gubler}) and prove a Poincar\'e lemma for these forms (Theorem \ref{PLemmaComplex}). The usual pullback of differential forms applied to superforms commutes with the differential operator along affine maps. For arbitrary maps this is false (cf. \ref{bem:counterexample}). Our proof of the Poincar\'e lemma follows the proof in the classical case. A crucial new tool is the introduction of a pullback of superforms via $C^\infty$-maps, which commutes with the differential operator. We use this to prove a homotopy formula in Theorem \ref{homotopyformula}. This will be the key result to prove the $d'$-Poincar\'e lemma. In section \ref{secfindim} we prove finiteness of the de Rham cohomology defined by superforms, using good covers and the Mayer-Vietoris sequence. In section \ref{secPLemmaComplex} we introduce real-valued differential forms on analytifications of algebraic varieties following Gubler's presentation in \cite{Gubler}. Then we use our result for polyhedral complexes to prove the Poincar\'e lemma for real-valued differential forms on the analytification of an algebraic variety (Theorem \ref{PLemmaXan}). Afterwards we sketch the argument for a generalization to a paracompact good analytic space (Theorem \ref{PLemmaCLD}). In Theorem \ref{ThmSkeleton} we show that, as a consequence, for a variety $X$ the cohomology of the complex $(A^{\bullet,0}(\Xan), d')$ is finite dimensional. \\
The author would like to thank to Walter Gubler, Johann Haas, Klaus K\"unnemann and Philipp Vollmer for reading various drafts of this work and providing very useful advise and the anonymous referee for his very precise review and helpful suggestions. The author would also like to thank the collaborative research centre SFB 1085 "Higher Invariants" by the Deutsche Forschungsgemeinschaft for its support. 

\newpage 
\section{A Poincar\'e lemma for superforms on polyhedral complexes} \label{secPLemma}
\subsection{Superforms on polyhedral complexes} \label{superformsoncomplexes}
Superforms were introduced by Lagerberg in \cite{Lagerberg} for open subsets of real vector spaces. They are analogues of $(p,q)$-forms on complex manifolds. The definition was extended to polyhedral complexes in \cite{CLD} (see also \cite{Gubler}). We recall the definitions.

\begin{defn} \label{defnd}
\begin{enumerate} [itemindent =*, leftmargin=0mm]
\item
For an open subset $U \subset \R^r$ denote by $A^{p}(U)$ the space of smooth real differential forms of degree $p$. Then the space of \textit{superforms of bidegree} $(p,q)$ on $U$ is defined as
\begin{align*}
A^{p,q}(U) := A^p(U) \otimes_{C^{\infty}(U)} A^q(U) = A^p(U) \otimes_\R \Lambda^q {\R^r}^*.
\end{align*}
If we choose a basis $x_1,\dots,x_r$ of $\R^r$ we can formally write a superform $\alpha \in A^{p,q}(U)$ as 
\begin{align*}
\alpha = \sum \limits _{|I| = p, |J| = q} \alpha_{IJ} d'x_I \vedge d''x_J
\end{align*}
where $I = \{ i_1 , \dots i_p\} $ and $J = \{j_1,\dots j_q\}$ are ordered subsets of $\{1,\dots,r\}$, $\alpha_{IJ} \in C^\infty(U)$ are smooth functions and 
\begin{align*}
d'x_I \vedge d''x_J := (dx_{i_1} \vedge \dots \vedge dx_{i_p}) \otimes_\R (dx_{j_1} \vedge \dots \vedge dx_{j_q}). 
\end{align*}
\item
There is a \textit{differential operator} $d' : A^{p,q} (U) = A^{p}(U) \otimes_\R \Lambda^q {\R^r}^* \rightarrow A^{p+1}(U) \otimes_\R \Lambda^q {\R^r}^* = A^{p+1,q}(U)$ which is given by $D \otimes \id$, where $D$ is the usual exterior derivative. We also have $A^{p,q}(U) = \Lambda^p {\R^r}^* \otimes_\R A^q(U)$ and can take the derivative in the second component. We put a sign on this operator and define $d'':= (-1)^p \id \otimes D$. In coordinates we have
\begin{align*}
d' \left( \sum \limits_{IJ} \alpha_{IJ} d'x_I \vedge d''x_J \right) = \sum \limits_{IJ} \sum \limits _{i = 1} ^{r} \frac {\del \alpha_{IJ} } {\del x_i} d'x_i \vedge d'x_I \vedge d''x_J
\end{align*} 
and 
\begin{align*}
d'' \left( \sum \limits_{IJ} \alpha_{IJ} d'x_I \vedge d''x_J \right) &= \sum \limits_{IJ} \sum \limits _{i = 1} ^{r} \frac {\del \alpha_{IJ} } {\del x_i} d''x_i \vedge d'x_I \vedge d''x_J \\
&= (-1)^p \sum \limits_{IJ} \sum \limits _{i = 1} ^{r} \frac {\del \alpha_{IJ} } {\del x_i} d'x_I \vedge  d''x_i \vedge d''x_J.
\end{align*} 
We further define $d := d' + d''$. The sign in $d''$ is such that $d'$ and $d''$ anticommute and hence $d$ is a differential. \\
\end{enumerate}
\end{defn}

\begin{bem} \label{bemsymmetry}
\begin{enumerate} [itemindent=*, leftmargin=0mm]
\item
There is an obvious symmetry in the definition of $d'$ and $d''$. If we switch factors in $A^{p,q} = A^{p} \otimes A^{q}$ then we change one into the other (up to sign). We will only talk about $d'$ in the following but corresponding statements are always true for $d''$. 
\item
The operator $d'$ is a differential. Hence for each $q \in \{0,\dots,r\}$ we get a complex
\begin{align*}
0 \rightarrow A^{0,q} \overset {d'} {\rightarrow} A^{1,q} \overset {d'} {\rightarrow} \dots \overset {d'} {\rightarrow} A^{r,q} \rightarrow 0
\end{align*}
of sheaves on $\R^r$.
\end{enumerate}
\end{bem}
 
\begin{bem} \label{defnaffinepullback}
If $F : \R^{r'} \rightarrow \R^{r}$ is an affine map and $U' \subset \R^{r}$ and $U \subset \R^{r}$ are open subsets such that $F(U') \subset U$, then there is a well defined \textit{pullback morphism} $F^*:A^{p,q}(U) \rightarrow A^{p,q}(U')$ that commutes with $d'$, $d''$ and $d$. 
\end{bem}

We will later define a pullback for a more general situation and use this in our proof of the Poincar\'e lemma. \\
Now we recall the definition of polyhedral complexes and forms on open subsets of polyhedral complexes following \cite{Gubler}. We refer to \cite[Appendix A]{Gubler2} for notations and results we use from convex geometry. 

\begin{defn} 
\begin{enumerate} [itemindent =*, leftmargin=0mm]
\item
A \textit{polyhedral complex} $\CC$ in $\R^r$ is a finite set of polyhedra (which we will always assume to be convex) in $\R^r$ with the following two properties: 
\begin{enumerate} [itemindent = \parindent]
\item
For  a polyhedron $\sigma \in \CC$, if $\tau$ is a face of $\sigma$ we have $\tau \in \CC$. 
\item
For two polyhedra $\sigma, \tau \in \CC$ we have that $\sigma \cap \tau$ is a face of both. 
\end{enumerate}
\item
The \textit{support} $|\CC|$ of $\CC$ is the union of all polyhedra in $\CC$. 
\item
A polyhedron $\sigma$ spans an affine space $\A_\sigma$ and we denote by $\Linear_\sigma$ the corresponding linear subspace of $\R^r$.
\item
Let $\Omega$ be an open subset of $|\CC|$. Then a \textit{superform $\alpha \in A^{p,q}(\Omega)$ of bidegree $(p,q)$} on $\Omega$ is given by a superform $\alpha' \in A^{p,q}(V)$ where $V$ is an open subset of $\R^r$ with $V \cap |\CC| = \Omega$. Two forms $\alpha' \in A^{p,q}(V)$ and $\alpha'' \in A^{p,q}(W)$ (with $V \cap |\CC| = W \cap |\CC| = \Omega)$ define the same form in $A^{p,q}(\Omega)$ if for each $\sigma \in \CC$ the restrictions of $\alpha'$ and $\alpha''$ to $\sigma \cap V = \sigma \cap W = \sigma \cap \Omega$ agree, which means that for all $x \in \sigma \cap \Omega$ and all tangent vectors $v_1,\dots,v_p,w_1,\dots,w_q \in \Linear_\sigma$ we have 
\begin{align*}
\langle \alpha'(x); v_1,\dots,v_p,w_1,\dots,w_q \rangle =  \langle \alpha''(x); v_1,\dots,v_p,w_1,\dots,w_q \rangle.
\end{align*}  
If  $\alpha \in A^{p,q}(\Omega)$ is given by $\alpha' \in A^{p,q}(V)$ we write $\alpha'|_\Omega = \alpha$. To simplify the notation we will often write $\alpha|_\sigma$ for $\alpha|_{\sigma \cap \Omega}$.
\end{enumerate}
\end{defn}

\begin{bem} \label{bemcomplex}
\begin{enumerate} [itemindent =*, leftmargin=0mm]
\item
We do not put any rationality assumption on our polyhedra, since the aspects of the theory of superforms we consider do not need it. These assumptions are needed in \cite{Gubler} to define integration of superforms, but we will not use this. 
\item
The set of superforms on an open subset $\Omega$ of a polyhedral complex $\CC$ depends only on the support of $\CC$. By this we mean that if $\DD$ is another polyhedral complex such that $\Omega$ is an open subset of $|\DD|$, then $A^{p,q}(\Omega)$ is the same whether we regard $\Omega$ as an open subset of $|\CC|$ or $|\DD|$.  
\item
The polyhedra in $\CC$ are partially ordered by the relation
\begin{align*}
\tau \prec \sigma :\Leftrightarrow \tau \text{ is a face of } \sigma.
\end{align*}
We will always assume our polyhedral complex to be of \textit{dimension n}, meaning that the maximal dimension of its polyhedra is $n$. In this case we have $A^{p,q}(\Omega) = 0$ for $\max(p,q) > n$. We say a polyhedral complex is \textit{pure of dimension n} if all maximal polyhedra are of dimension $n$.
\item
Taking $d'$ of a superform on an open subset of $\R^r$ is compatible with restriction to polyhedra. Hence the differential $d'$ induces, for an open subset $\Omega \subset |\CC|$, a differential $d': A^{p,q}(\Omega) \rightarrow A^{p+1,q}(\Omega)$. 
\item
A partition of unity argument shows that $A^{p,q}$ is indeed a sheaf on $|\CC|$ and hence for each $q \in \{0,\dots,n\}$ we get a complex
\begin{align*}
0 \rightarrow A^{0,q} \overset {d'} {\rightarrow} A^{1,q} \overset {d'} {\rightarrow} \dots \overset {d'} {\rightarrow} A^{n,q} \rightarrow 0
\end{align*}
of sheaves on $|\CC|$. The fact that this complex is exact in positive degrees will be the main result of this section. 
\item
We define $A^{k}(\Omega) := \bigoplus \limits_{p+q = k} A^{p,q}(\Omega)$. Thus $d: A^k(\Omega) \rightarrow A^{k+1}(\Omega)$ is a differential. We get a complex
\begin{align*}
0 \rightarrow A^0 \overset {d} {\rightarrow} A^1 \overset {d}{\rightarrow} \dots \overset {d} {\rightarrow} A^{2n} \rightarrow 0
\end{align*}
of sheaves on $|\CC|$. 
\item
The affine pullback as in Remark \ref{defnaffinepullback} is compatible with restriction to polyhedra. Hence if $F: \R^{r'} \rightarrow \R^r$ is an affine map, $\CC$ resp. $\CC'$ are polyhedral complexes in $\R^r$ resp. $\R^{r'}$  such that $F(|\CC'|) \subset |\CC|$ and $\Omega \subset |\CC|$ resp. $\Omega' \subset |\CC'|$ are open subsets such that $F(\Omega') \subset \Omega$ then the affine pullback induces a well defined pullback $F^*: A^{p,q}(\Omega) \rightarrow A^{p,q}(\Omega')$.
\end{enumerate}
\end{bem}

\subsection {A $d'$-Poincar\'e lemma for superforms on polyhedral complexes}

In this subsection we will prove a $d'$-Poincar\'e lemma for superforms on polyhedral complexes. The polyhedral complex $\CC$ will always be of dimension $n$.

\begin{lem} [Chain Homotopy Lemma]  \label{d'Homlemma}
Let $\CC$ be a polyhedral complex in $\R^{r}$ and $\Omega \subset |\CC|$ an open subset. Let $B = [0,1] \subset \R$ be the closed unit interval and for $i =0,1$ 
\begin{align*}
\iota_{i}: \Omega \rightarrow \Omega \times \{i\} \subset \Omega \times B
\end{align*}
the inclusions. Then for all $p \in \{0,\dots,n+1\}$ and $q \in \{0,\dots,n\}$ there exists a linear map
\begin{align} \label{defnK}
 K' : A^{p,q} (\Omega  \times B) \rightarrow A^{p-1, q} (\Omega),
\end{align}
such that
\begin{align} \label{equationforK}
d'K' + K'd' = \iota_1^* - \iota_0^*.
\end{align}
\begin{proof}
The proof is a variant of the classical chain homotopy lemma for ordinary differential forms. Observe first that $|\CC| \times B$ is the support of a polyhedral complex in $\R^r \times \R$ and hence it makes sense to talk about superforms on $\Omega \times B$. Let $\alpha \in A^{p,q}(\Omega \times B)$ be given by $\beta \in A^{p,q} (V \times B')$ for some open set $V \subset \R^{r}$ and some open interval $B'$ such that $B \subset B' \subset \R$. Let $x_1,\dots,x_r$ be a basis of $\R^r$ and denote by $t$ the coordinate of $B$. We write 
\begin{align} \label{decomposition}
\beta = &\sum \limits _{|I| = p, |J| = q } a_{IJ} d'x_I \vedge d''x_J \\
\notag&+ \sum \limits_{|I| = p-1, |J| = q} b_{IJ} d't \vedge d'x_I \vedge d''x_J  \\
\notag&+  \sum \limits_{|I| = p, |J| = q-1} e_{IJ} d'x_I \vedge d''t \vedge d''x_J \\
\notag&+  \sum \limits_{|I| = p-1, |J| = q-1 } g_{IJ} d't \vedge d'x_I \vedge d''t \vedge d''x_J.
\end{align}
Then we define 
\begin{align*}
K': A^{p,q}(V \times B) &\rightarrow A^{p - 1,q}(V) \\
\beta &\mapsto  \sum \limits _{|I| = p-1, |J| = q} c_{IJ} d'x_I \vedge d''x_J \\
\text{ with } c_{IJ}(x) :&= \int \limits _0 ^1 b_{IJ} (x, t) dt.
\end{align*}
We show that this definition is independent of the choice of the basis $x_1,\dots,x_r$. Let therefore $y_1,\dots,y_r$ be another basis. First of all we notice that the decomposition into the four summands as in (\ref{decomposition}) is not affected by our base change. We further notice that 
\begin{align*}
d'x_I \vedge d''x_J = \sum \limits_{|I'| = |I| ,|J'| = |J|} \lambda_{I,I'} \lambda_{J,J'} d'y_{I'} \vedge d''y_{J'},
\end{align*}
where $\lambda_{I,I'}$ is the determinant of the $I \times I'$ minor of the base change matrix from $x_1,\dots,x_r$ to $y_1,\dots,y_r$ and similar for $J$ and $J'$. Now we have
\begin{align*}
b_{IJ}  d't \vedge d'x_I \vedge d''x_J = b_{IJ} \sum \limits_{I',J'} \lambda_{I,I'} \lambda_{J,J'} d't \vedge d'y_{I'} \vedge d''y_{J'}
\end{align*}
and this term is mapped under $K'$ to
\begin{align*}
&\sum \limits_{I',J'} \left(\int \limits _0 ^1   \lambda_{I,I'} \lambda_{J,J'} b_{IJ} dt \right) d'y_{I'} \vedge d''y_{J'} \\
=  &\left(\int \limits _0 ^1   b_{IJ} dt \right) \sum \limits_{I',J'}\lambda_{I,I'} \lambda_{J,J'} d'y_{I'} \vedge d''y_{J'} \\
=  &\left(\int \limits _0 ^1   b_{IJ} dt \right) d'x_I \vedge d''x_J,
\end{align*}
which shows the independence on the choice of the basis. \\
Given $V$ and $B'$ we have the diagram
\begin{align*}
\begin{xy}
\xymatrix 
{ 
A^{p,q}(V \times B') \ar[rr]^{K'} \ar[d] && A^{p-1, q} (V) \ar[d] \\ 
A^{p,q} (\Omega \times B) \ar@{-->}[rr]&& A^{p-1,q} (\Omega).
}
\end{xy}
\end{align*}
To get a well defined map on the bottom that makes this diagram commutative, we need that $\beta|_{\sigma \times B} = 0$ for all $\sigma \in \CC$ implies $K'(\beta)|_{\sigma} = 0$ for all $\sigma \in \CC$. Let therefore $\sigma$ be a maximal polyhedron in $\CC$ and $W = V \cap \sigma$. It suffices to show that if $\beta|_{W \times B} = 0$, then $K'(\beta)|_W = 0$. By what we did above we may choose a basis as we like. Let therefore $x_1,\dots,x_m$ be a basis of $\Linear_\sigma$ and $x_{m+1},\dots,x_r$ a basis of a complement. Then from $\beta |_{W \times B} = 0$ we get $b_{IJ}|_{W \times B} = 0$ for all $I, J \subset \{1,\dots,m\}$. This means however that $c_{IJ}|_W = 0$ for all $I,J \subset \{1,\dots,m \}$. From that we get $K'(\beta)|_W = 0$. Hence setting $K'(\alpha) := K'(\beta)$ is independent of the choice of the form $\beta$ by which $\alpha$ is given. It is also independent of the choice of $V$ and $B'$. This gives a well defined map
\begin{align*}
K': A^{p,q} (\Omega \times B) \rightarrow A^{p-1,q} (\Omega)
\end{align*}
as required in (\ref{defnK}). We will now show that (\ref{equationforK}) holds. It is enough to check that 
\begin{align*}
d'K' \beta + K' d' \beta = \iota_1^*\beta - \iota_0^*\beta
\end{align*}
holds for every $\beta \in A^{p,q}(V \times B')$, where $V$ is an open subset of $\R^{r}$ and $B'$ is an open interval such that $B \subset B' \subset \R$. 
It suffices to check the following four cases:
\begin{enumerate} [itemindent =*, leftmargin=0mm]
\item 
$\beta = a_{IJ} d'x_I \vedge d''x_J$: \\
We have $K'(\beta) = 0$ and 
\begin{align*}
K'(d'(\beta)) &= K'  \left(\frac {\del a_{IJ}} {\del t} d't \vedge d'x_I \vedge d''x_J \right) \\
&+ \sum \limits_{i=1} ^{r} K' \left( \frac {\del a_{IJ}} {\del x_i} d'x_i \vedge d'x_I \vedge d''x_J \right) \\
&= \left(\int \limits _0 ^1 \frac {\del a_{IJ}} {\del t}dt \right) d'x_I \vedge d''x_J \\
&= (a_{IJ}(. ,1) - a_{IJ} (. ,0))  d'x_I \vedge d''x_J \\
&= \iota_1^*(\beta) - \iota_0^*(\beta)
\end{align*}
\item
$\beta = b_{IJ} d't \vedge d'x_I \vedge d''x_J$: \\
We have $\iota_1^* (\beta) = \iota_0^* (\beta) = 0$, since the pullback of $d't$ is zero. We further have
\begin{align*}
d'K'(\beta) = \sum \limits_{i = 1} ^{r} \left(\int \limits ^1 _0 \frac {\del b_{IJ}} {\del x_i}dt \right) d'x_i \vedge d'x_I \vedge d''x_J 
\end{align*}
and
\begin{align*}
K'd'(\beta) &= \sum \limits _{i= 1} ^{r} K' \left( \frac {\del b_{IJ}} {\del x_i} d'x_i \vedge d't \vedge d'x_I \vedge d''x_J \right)\\
&= - \sum \limits _{i= 1} ^{r} K' \left(\frac {\del b_{IJ}} {\del x_i} d't \vedge d'x_i \vedge d'x_I \vedge d''x_J\right) \\
&= - \sum \limits _{i= 1} ^{r}  \left(\int   \limits ^1 _0 \frac {\del b_{IJ}} {\del x_i} dt \right)  d'x_i  \vedge d'x_I \vedge d''x_J.
\end{align*}
\item
$\beta =  e_{IJ} d'x_I \vedge d''t \vedge d''x_J$ :\\
Similarly to ii) the pullbacks are zero. Since both $\beta$ and $d'\beta$ have a factor $d''t$ by definition they are sent to $0$ by $K'$.
\item
$\beta =  g_{IJ} d't \vedge d'x_I \vedge d''t \vedge d''x_J$: \\
Same as iii).
\end{enumerate}
Adding up these parts we have proven that (\ref{equationforK}) holds on $V$. Now if $\alpha \in A^{p,q}(\Omega \times B)$ is given by $\beta \in A^{p,q}(V \times B')$ then the equation holds for $\alpha$ simply because it holds for $\beta$.
\end{proof}
\end{lem}

In the classical proof of the Poincar\'e lemma for star shaped subsets $U$ of $\R^n$ the idea is to pull back differential forms via a contraction of $U$ to its center. This contraction is however not an affine map. So we will introduce in Definition \ref{smoothpullback} a pullback for superforms along $C^{\infty}$-maps that still commutes with $d'$ (as we will see in \ref{bemsmoothpullback}). This will be a crucial ingredient in our proof of the Poincar\'e lemma for superforms. The following example shows that the direct approach does not work. 

\begin{bem} \label{bem:counterexample}
Given a $C^\infty$-map $F: V' \rightarrow V$, where $V'$ resp. $V$ are open subsets of $\R^{r'}$ resp. $\R^r$ we can define a naive pullback 
\begin{align*}
F^*: A^{p,q}(V) = A^{p}(V) \otimes A^{q}(V) \rightarrow A^p(V') \otimes A^q(V') = A^{p,q}(V'),
\end{align*}
which is just given by the tensor products of the usual pullback of differential forms. This pullback however does not commute with the differential $d'$ in general, as can be seen in the following example. Let $V' = \R^2$ and $V = \R$ and $F(x,y) = xy$. Denote the coordinate on $\R$ by $t$. Then we have $d' F^*(d''t) =d'(x d''y + y d''x) = d'x \vedge d''y + d'y \vedge d''x \neq 0$, however $d' (d''t) = 0$ and thus $F^*(d' d''t) = 0$.  \\
The reason for this is that the definition of this pullback uses the presentation $A^{p,q} = A^p \otimes A^q$, while the definition of $d'$ uses the presentation $A^{p,q} = A^{p} \otimes \Lambda^q {\R^r}^*$ and thus these two are not compatible. We would therefore like to define a pullback which uses the presentation $A^{p,q} = A^{p} \otimes \Lambda^q {\R^r}^*$. 
\end{bem}

\begin{lem} \label{lemressurj} 
Let $\CC$ be a polyhedral complex in $\R^{r}$, $\Omega \subset |\CC|$ an open subset and $W \subset \R^r$ an open subset such that $\Omega = W \cap |\CC|$. Then the restriction map $A^{p,q}(W) \rightarrow A^{p,q}(\Omega)$ is surjective. In particular, we may assume that any form $\alpha \in A^{p,q}(\Omega)$ is given by a form on $W$.  
\begin{proof}
Let $\alpha \in A^{p,q}(\Omega)$ be given by $\beta \in A^{p,q}(V)$. Then $\alpha$ is also given by $\beta|_{V \cap W}$, hence we may assume $V \subset W$. Notice that $\Omega$ is a closed subset of $W$. Choose a function $f \in C^\infty(W)$ such that $f|_\Omega \equiv 1$ and $\supp_W f \subset V$. Then $\alpha$ is given by $f|_V \beta$ and this can be extended by zero to a form in $A^{p,q}(W)$. 
\end{proof}
\end{lem}

\begin{defn} [$C^\infty$-pullback of $(p,q)$-forms] \label{smoothpullback} 
We define a pullback for superforms on open subsets $V \subset \R^{r}$ and, under certain conditions, for superforms on polyhedral complexes. 
\begin{enumerate} [itemindent =*, leftmargin=0mm]
\item
Let $V' \subset \R^{r'}$ and $V \subset \R^{r}$ be open subsets. Let $F = (s_F, L_F)$ be a pair of maps such that $s_F:V' \rightarrow V$ is a $C^{\infty}$-map and $L_F: \R^{r'} \rightarrow \R^{r}$ is linear. We define
\begin{align*}
F^* := s_F^* \otimes L_F^* : A^{p,q} (V) = &A^p(V) \otimes_\R \Lambda^q {\R^r}^* \\
\rightarrow &A^{p}(V') \otimes_\R \Lambda^q {\R^{r'}}^* = A^{p,q}(V').
\end{align*}
Explicitly, if $\beta \in A^{p,q}(V)$ we have
\begin{align*}
&\langle F^*(\beta)(x); v_1,\dots,v_p,w_1,\dots,w_q \rangle \\
 = &\langle \beta(s_F(x)); d(s_F)_x(v_1),\dots,d(s_F)_x(v_p),L_F(w_1),\dots,L_F(w_q) \rangle 
\end{align*}
for all $x \in V'$ and $v_i, w_i \in \R^{r'}$, where $d(s_F)_x$ denotes the differential of $s_F$ at $x$. 
\item 
Let $\CC'$ and $\CC$ be polyhedral complexes in $\R^{r'}$ and $\R^{r}$ respectively. Let $\Omega' \subset |\CC'|$ and $\Omega \subset |\CC|$ open subsets and $V'$ resp. $V$ open neighbourhoods of $\Omega'$ resp. $\Omega$ in $\R^{r'}$ resp. $\R^r$. Let $s_F: V' \rightarrow V$ be a $C^\infty$-map and $L_F: \R^{r'} \rightarrow \R^r$ a linear map such that $s_F(\Omega') \subset \Omega$.  The pair $F = (s_F, L_F)$ is said to \textit{allow a pullback from $\Omega$ to $\Omega'$} if there exist open subsets $W$ of $V$ and $W'$ of $V'$ such that $W \cap |\CC| = \Omega$, $W' \cap |\CC'| = \Omega'$, $s_F(W') \subset W$ and for all $\beta \in A^{p,q}(W)$ such that $\beta|_\Omega = 0$ we have $F^*(\beta)|_{\Omega'} = 0$. In that case, for a form $\alpha \in A^{p,q}(\Omega)$ we choose $\beta \in A^{p,q}(W)$ by which $\alpha$ is given (which is possible by \ref{lemressurj}) and we define $F^*(\alpha) \in A^{p,q}(\Omega')$ to be given by $F^*(\beta) \in A^{p,q}(W')$. The form $F^*(\alpha) \in A^{p,q}(\Omega')$ is then independent of the choice of $W$, $W'$ and $\beta$, as will be shown in the next Lemma.  
\end{enumerate}
\end{defn}

\begin{lem}
The definition of $F^*(\alpha)$ above is independent of the choice of $W, W'$ and $\beta$. 
\begin{proof}
The independence of $\beta$ is simply due to the property that $F^*$ respects forms that restrict to zero. \\
Now if both $W_1, W'_1, \beta_1$ and $W_2,W'_2, \beta_2$ have the properties required in the definition above, then by \ref{lemressurj} we can choose a form $\delta \in A^{p,q}(W_1  \cup W_2)$ such that $\delta|_{\Omega} = \alpha$. By independence of the form we have
\begin{align*}
F^*(\beta_1)|_{\Omega'} = F^*(\delta|_{W_1})|_{\Omega'} = F^*(\delta)|_{W'_1}|_{\Omega'} =  F^*(\delta)|_{\Omega'}
\end{align*}
and the same works for $F^*(\beta_2)|_{\Omega'}$, which proves exactly the independence we wanted to show.
\end{proof}
\end{lem}

\begin{bem} \label{bemsmoothpullback} 
\begin{enumerate} [itemindent =*, leftmargin=0mm]
\item
The pullback between open subsets of vector spaces commutes with taking $d'$ since both use the presentation $A^{p,q}(V) = A^{p}(V) \otimes {\R^r}^*$. We have $F^* = s_F^* \otimes L_F^*$ and $d' = D \otimes \id$ and $s_F^*$ and $D$ commute. If $F$ allows a pullback, then the pullback $F^*$ between open subsets of the supports of polyhedral complexes commutes with $d'$ since both $F^*$ and $d''$ are defined via restriction. 
\item 
The pullback is functorial in the following sense: Let $\CC$, $\CC'$ and $\CC''$ be polyhedral complexes, $\Omega \subset |\CC|$, $\Omega' \subset |\CC'|$ and $\Omega'' \subset |\CC''|$ open subsets and $V \subset \R^{r}$ resp. $V' \subset \R^{r'}$ resp. $V'' \subset \R^{r''}$ open neighbourhoods of $\Omega$ resp $\Omega'$ resp. $\Omega''$. Let further $F = (s_F, L_F)$ and $G = (s_G, L_G)$ be pairs of maps such that $s_F: V' \rightarrow V$ and $s_G: V'' \rightarrow V'$ are $C^\infty$-maps, $L_F: \R^{r'} \rightarrow \R^{r}$ and $L_G: \R^{r''} \rightarrow \R^{r'}$ are linear maps and $s_F (\Omega') \subset \Omega'$ and $s_G(\Omega'') \subset \Omega'$. If both $F$ resp. $G$ allow a pullback from $\Omega$ to $\Omega'$ resp. $\Omega'$ to $\Omega''$ and we define $F \circ G := (s_F \circ s_G, L_F \circ L_G)$ then $F \circ G$ allows a pullback from $\Omega$ to $\Omega''$ and we have $(F \circ G)^* = G^* \circ F^*$.  
\item
Let $F: \R^{r'} \rightarrow \R^r$ be an affine map and denote by $\Linear_F := F - F(0)$ the associated linear map. Then the pullback via $F$ in the sense of Remark \ref{defnaffinepullback} is the pullback via $(F , \Linear_F)$ in the sense of Definition \ref{smoothpullback} above. 
\item
The notion of allowing pullback does not depend on the underlying polyhedral complexes $\CC$ and $\CC'$, since the pullback is defined on an open neighborhood and the restriction only depends on $\Omega$ and $\Omega'$ (cf. \ref{bemcomplex}).   
\end{enumerate}
\end{bem}

\begin{satz} [Homotopy Formula] \label{homotopyformula}
Let $V$ an open subset of $\R^r$. Let further $s_F: V \rightarrow V$ be a $C^\infty$-map and $L_F := \id$. Let $s_G : V \times \R \rightarrow V$ such that $s_G(.,0) = s_F$ and $s_G(.,1) = \id$. Let $L_G = \pr_1 : \R^r \times \R \rightarrow \R^r$ be the projection to the first factor. Denote by $F^*$ respectively $G^*$ the pullback from $V$ to $V$ respectively to $V \times \R$ via pairs $(s_F, L_F)$ respectively $(s_G, L_G)$. Then for $\alpha \in A^{p,q}(V)$ we have
\begin{align*}
\alpha - F^* \alpha  = d' K' G^* \alpha + K' G^* d' \alpha
\end{align*}
for any operator $K'$ satisfying the equality (\ref{equationforK}) of Lemma \ref{d'Homlemma}. 
\begin{proof}
We calculate 
\begin{align*}
\id^* - F^*  &= (G \circ \iota_1)^* - (G \circ \iota_0)^* \\
&= \iota_1^* \circ G^*- \iota_0^* \circ G^* \\
&= (\iota_1^* - \iota_0^*)  \circ G^* \\
&= (K'd' + d'K') G^*\\
&= K'd' G^* + d' K' G^* \\
& \overset {(\ref{bemsmoothpullback})} {=} K' G^* d' + d' K' G^* ,
\end{align*}
where we denote by $\iota_i$ the pair $(\iota_i, \Linear_{\iota_i})$. Now putting in $\alpha$ and using $\id^*(\alpha) = \alpha$ gives the desired result.
\end{proof}
\end{satz}

Note that if $\Omega$ is an open subset of $|\CC|$ for some polyhedral complex $\CC$ in $\R^r$ and $F$ resp. $G$ allows a pullback from $\Omega$ to $\Omega$ resp. $\Omega \times B$, where $B = [0,1]$ is the closed unit interval, then the analogue formula also holds for $\alpha \in A^{p,q}(\Omega)$ since all operators are defined via restriction.

\begin{defn}
Let $\CC$ be a polyhedral complex in $\R^r$. An open subset $\Omega$ of $|\CC|$ is called \textit{polyhedrally star shaped} with centre $z$ if there is a polyhedral complex $\DD$ such that $\Omega$ is an open subset of $|\DD|$ and for all maximal $\sigma \in \DD$ the set $ \sigma \cap \Omega$ is star shaped with centre $z$ in the sense that for all $x \in \sigma \cap \Omega$ and for all $t \in [0,1]$ the point $z + t(x-z)$ is contained in $\sigma \cap \Omega$ . 
\end{defn}

\begin{bem}
It is obvious that if $\Omega \subset \CC$ is a polyhedrally star shaped open subset with centre $z$, then $\Omega$ is also star shaped with centre $z$. The converse is not true however: Take $\CC$ such that $|\CC| = [-1,1] \times [-1,1] \cup \{0\} \times  [1,2] \cup [1,2] \times \{0\} \subset \R^2$. Then $\Omega := |\CC|$ is star shaped but not polyhedrally star shaped.
\end{bem}

\begin{lem} \label{lemallowpullback}
Let $\CC'$ and $\CC$ be polyhedral complexes in $\R^{r'}$ and $\R^r$ respectively. Let $\Omega' \subset |\CC'|$ and $\Omega \subset |\CC|$ be open subsets and $V'$ resp. $V$ open neighbourhoods of $\Omega'$ resp. $\Omega$ in $\R^{r'}$ resp. $\R^r$. Let $s_F: V' \rightarrow V$ be a $C^\infty$-map and $L_F: \R^{r'} \rightarrow \R^r$ a linear map such that $s_F(\Omega') \subset \Omega$. Suppose there exist polyhedral complexes $\DD'$ in $\R^{r'}$ and $\DD$  in $\R^{r}$ such that $\Omega'$ resp. $\Omega$ are open subsets of $|\DD'|$ resp. $|\DD|$ and such that for all maximal $\sigma' \in \DD'$ there exists a maximal $\sigma \in \DD$ such that we have
\begin{enumerate} [itemindent =*, leftmargin=0mm]
\item [(a)]
$\forall x \in \sigma' \cap \Omega', s_F(x) \in \sigma$ and 
\item [(b)]
$ \forall w \in \Linear_{\sigma'}, L_F(w) \in \Linear_\sigma$.
\end{enumerate}
Then $F := (s_F, L_F)$ allows a pullback from $\Omega$ to $\Omega'$.
\begin{proof}
We  first note that Remark \ref{bemsmoothpullback} says that whether that $F$ allows a pullback does not depend on the polyhedral complex and hence we may assume that $\DD = \CC$ and $\DD' = \CC'$. Let $W \subset V$ be an open subset such that $W \cap |\CC| = \Omega$ and let $\beta \in A^{p,q}(W)$. For $F$ to allow a pullback we have to show that if $\beta |_\sigma = 0$ for all maximal polyhedra $\sigma \in \CC$ then $(F^*\beta)|_{\sigma'} = 0$ for all maximal $\sigma' \in \CC'$. \\
Let $\sigma' \in \CC'$ be a maximal polyhedron and $\sigma \in \CC$ the maximal polyhedron such that $\sigma$ and $\sigma'$ satisfy condictions (a) and (b). We then have that 
\begin{enumerate} [itemindent =*, leftmargin=0mm]
\item [(c)]
$ \forall x \in \sigma' \cap \Omega', \; \; \forall v \in \Linear_{\sigma'}, d(s_F)_x(v) \in \Linear_\sigma$,
\end{enumerate}
due to condition (a).\\
For $\sigma \in \CC$ the fact that $\beta|_{\sigma} = 0$ just means 
\begin{align*}
\langle \beta(x); v_1,\dots ,v_p,w_1,\dots,w_q \rangle = 0
\end{align*}
 for all $x \in \sigma \cap \Omega , v_i, w_i \in \Linear_{\sigma}$. But then we have  
\begin{align*}
&\langle F^*(\beta)(x); v_1,\dots,v_p,w_1,\dots,w_q \rangle \\
 = &\langle \beta(s_F(x)); d(s_F)_x(v_1),\dots,d(s_F)_x(v_p),L_F(w_1),\dots,L_F(w_q) \rangle = 0
\end{align*}
for all $x \in \sigma' \cap\Omega', v_i, w_i \in \Linear_{\sigma'}$ by conditions (a), (b) and (c). Hence $F^*(\beta)|_{\sigma'} = 0$. This shows that if $\beta|_\sigma = 0$ for all $\sigma \in \CC$ then $F^*(\beta)|_{\sigma'} = 0$ for all maximal and hence all $\sigma' \in \CC'$. Thus $F$ allows a pullback from $\Omega$ to $\Omega'$. 
\end{proof}
\end{lem}

\begin{satz} [$d'$-Poincar\'e lemma for polyhedral complexes]   \label{PLemmaComplex}
Let $\CC$  be a polyhedral complex in $\R^r$ and $\Omega \subset |\CC|$ a polyhedrally star shaped open subset with centre $z$. Let $\alpha \in A^{p,q}(\Omega)$ with $p> 0$ and $d' \alpha = 0$. Then there exists $\beta \in A^{p-1,q}(\Omega)$ such that $d'\beta = \alpha$. 
\begin{proof}
We want to use Theorem \ref{homotopyformula} with $s_F$ the constant map to the centre $z$ of $\Omega$. Let $L_F = \id$, $s_G$ given by
\begin{align*}
s_G : \R^r \times \R &\rightarrow \R^r \\
(x , t) &\mapsto z + t (x - z)
\end{align*}
and  $L_G =  \pr_1$. It is easy to check that both $F$ and $G$ have the properties required in Theorem \ref{homotopyformula}. We show that they allow a pullback from $\Omega$ to $\Omega$ resp. $\Omega \times B$ by showing that they fulfill the conditions required in Lemma \ref{lemallowpullback}. Since $\Omega$ is polyhedrally star shaped we know that there exists a polyhedral complex $\DD$ such that $\Omega$ is an open subset of $|\DD|$ and such that $\sigma \cap \Omega$ is star shaped with centre $z$ for all maximal $\sigma \in \DD$. We take $\DD'$ to be the polyhedral complex whose maximal polyhedra are of the form $\sigma \times B$ for $\sigma \in \DD$ a maximal polyhedron. Let $\sigma' = \sigma \times B \in \DD'$ be such a maximal polyhedron. For $(x,t) \in \sigma'$ we have $s_G(x,t) \in \sigma$ because $\sigma \cap \Omega$ is star shaped with centre $z$. Since it is obvious that $L_G(\Linear_{\sigma'}) \subset \Linear_\sigma$, $G$ allows a pullback from $\Omega$ to $\Omega \times B$ by Lemma \ref{lemallowpullback}. Since $s_F$ is constant and $L_F$ is the identity we also see that $F$ has the properties of Lemma \ref{lemallowpullback} and hence allows a pullback from $\Omega$ to $\Omega$.
Now since $\alpha \in A^{p,q}(\Omega)$ with $p > 0$ we have $F^*\alpha = 0$ (since $s_F$ is a constant map). Together with our assumption $d'\alpha = 0$, Theorem \ref{homotopyformula} yields
\begin{align*}
\alpha =  d'(K'G^*\alpha),
\end{align*}
which proves the theorem.
\end{proof}
\end{satz}

\begin{bem} \label{bemlocallystarshaped}
Let $\CC$ be a polyhedral complex, $\Omega \subset |\CC|$ an open subset and $z \in \Omega$ a point. Let $V \subset \R^r$ be an open ball around $z$ such that $V \cap |\CC| \subset \Omega$ and such that $V$ intersects only polyhedra in $\CC$ that contain $z$. Write $\Omega' := V \cap |\CC|$ and let $\DD$ be the polyhedral subcomplex of $\CC$ whose maximal polyhedra are the ones in $\CC$ which intersect $\Omega'$. Let $\sigma$ be a maximal polyhedron in $\DD$. Since $V$ and $\sigma$ are both convex, $\sigma \cap V = \sigma \cap \Omega'$ is convex and hence it is star shaped with respect to any point and in particular with respect to $z$. Hence any point $z \in |\CC|$ has a basis of open neighbourhoods consisting of polyhedrally star shaped open sets.
\end{bem}

\begin{kor} \label{korplemmacomplex}
For all $q \in \{0,\dots,n\}$ the complex 
\begin{align*}
0 \rightarrow A^{0,q} \overset {d'} {\rightarrow} A^{1,q} \overset {d'} {\rightarrow} \dots \overset {d'} {\rightarrow} A^{n,q} \rightarrow 0
\end{align*}
of sheaves on $|\CC|$ is exact in positive degrees. 
\begin{proof}
This is a direct consequence of Theorem \ref{PLemmaComplex} and Remark \ref{bemlocallystarshaped}.
\end{proof}
\end{kor}

\begin{kor} \label{korPLemmaComplexII}
The complex
\begin{align} \label{acycrescomplex}
0 \rightarrow \underline{\R} \rightarrow A^{0,0} \overset{d'} {\rightarrow} A^{1,0} \overset {d'} {\rightarrow} \dots \overset{d'} {\rightarrow}  A^{n,0} \rightarrow 0
\end{align}
of sheaves on $|\CC|$ is exact. The cohomology of its complex of global sections
\begin{align} \label{calcohocomplex}
0 \rightarrow A^{0,0}(|\CC|) \overset{d'} {\rightarrow} A^{1,0}(|\CC|) \overset {d'} {\rightarrow} \dots \overset{d'} {\rightarrow}  A^{n,0}(|\CC|) \rightarrow 0
\end{align}
is isomorphic to the sheaf cohomology $H^*(|\CC|, \underline{\R})$ of the constant sheaf $\underline{\R}$ and to the singular cohomology $H_{\sing}^*(|\CC|, \R)$.
\begin{proof}
The fact that $f \in A^{0,0}(\Omega)$ for $\Omega \subset |\CC|$ is a function and that $d'f = 0$ if and only if $f$ is locally constant together with Corollary $\ref{korplemmacomplex}$ show that the complex (\ref{acycrescomplex}) is exact. Since the sheaves $A^{p,q}$ admit partitions of unity, they are fine, hence acyclic \cite[Chapter II, Proposition 3.5 and Theorem 3.11]{Wells}. This means that the complex (\ref{calcohocomplex}) calculates the sheaf cohomology of $\underline{\R}$. Since $|\CC|$ is paracompact, Hausdorff and locally compact, this is the singular cohomology of the topological space $|\CC|$ \cite[Chapter III, Theorem 1.1]{Bredon}.
\end{proof}
\end{kor}

\begin{bem} \label{remnotexact}
We can not expect a $d$-Poincar\'e lemma to hold for the following reason: If $J: A^{p,q} \rightarrow A^{q,p}$ is the operator that switches the factors in the tensor product $A^{p,q} = A^{p} \otimes_{C^\infty} A^q$, then for any function $f \in A^{0,0}(V)$, where $V$ is an open subset of $\R^r$, we have that $df$ is invariant under $J$ but there is no need for a $d$-closed 1-form to be invariant under $J$. 
\end{bem}

\section{Finiteness results for the cohomology of superforms on polyhedral complexes} \label{secfindim}

\subsection{Good covers and Mayer-Vietoris-Sequence}

In this subsection $\CC$ is a polyhedral complex in $\R^r$. Recall from Remark \ref{bemcomplex} that $A^k = \bigoplus \limits_{p+q = k} A^{p,q}$. Let $(A^\bullet, D)$ denote either the complex $(A^\bullet,d)$ or the complex $(A^{\bullet,q}, d')$, for fixed $q$, of sheaves of superforms on $|\CC|$. For an open subset $\Omega \subset |\CC|$ we will write $H^{\bullet}(\Omega)$ for the cohomology of $(A^\bullet(\Omega), D)$. The statements of Theorem \ref{MVS} and Lemma \ref{lemgoodcover} are special cases of theorems which are certainly well known. We choose to present them here with short proofs for the convenience of the reader. 

\begin{defn} \label{defngoodcover}
Let $\Omega \subset |\CC|$ be an open subset. An open cover $(\Omega_i)_{i \in I}$  of $\Omega$ is called a \textit{reasonable cover for $(A^\bullet, D)$} if $I$ is finite and for all $n \in \N_{>0}$ and for all $\iota_1,\dots,\iota_n \in I$ the set $\Omega_{\iota_1,\dots,\iota_n} := \bigcap \limits_{i=1} ^n \Omega_{\iota_i}$ has the property that $H^{k} (\Omega_{\iota_1,\dots,\iota_n})$ is finite dimensional for all $k \in \N_0$. It is called a \textit{good cover} if further $H^{k}(\Omega_{\iota_1,\dots,\iota_n}) = 0$ for all $k >0$. 
\end{defn}

\begin{satz} [Mayer-Vietoris-Sequence] \label{MVS}
Let $\CC$ be a polyhedral complex and $\Omega, \Omega_1, \Omega_2$ open subsets of $|\CC|$ such that $\Omega = \Omega_1 \cup \Omega_2$. Let further $\Omega_{12} := \Omega_1 \cap \Omega_2$. Then there exists a long exact sequence 
\begin{align*}
0 \rightarrow &H^{0}(\Omega) \rightarrow H^{0}(\Omega_1) \oplus H^{0} (\Omega_2) \rightarrow H^{0}(\Omega_{12}) \rightarrow \dots \\
\dots H^{k-1}(\Omega_{12}) \rightarrow  &H^{k}(\Omega) \rightarrow H^{k}(\Omega_1) \oplus H^{k} (\Omega_2) \rightarrow H^{k}(\Omega_{12}) \rightarrow \\
 \rightarrow  &H^{k+1}(\Omega) \rightarrow \dots 
\end{align*}
\begin{proof}
A partition of unity argument shows that the sequence
\begin{align*}
0 \rightarrow A^\bullet(\Omega) \rightarrow A^\bullet(\Omega_1) \oplus A^\bullet(\Omega_2) \rightarrow A^\bullet(\Omega_{12}) \rightarrow 0
\end{align*}
is exact. The result is then obtained by the long exact cohomology sequence.
\end{proof}
\end{satz}

\begin{lem} \label{lemgoodcover}
Let $\Omega \subset |\CC|$ be an open subset and $(\Omega_i)_{i=1,\dots,m}$ be a reasonable cover of $\Omega$ for $A^\bullet$. Then $H^{k}(\Omega)$ is a finite dimensional real vector space for all $k \in \N_0$.\\
If  $(\Omega_i)_{i=1,\dots,m}$ is a good cover, then we further have $H^{k}(\Omega) = 0$ if $k \geq m$. 
\begin{proof}
We use induction on $m$ with $m=1$ being just the definition \ref{defngoodcover} in both the reasonable and the good case. Now let $m \geq 2$. Let $\Omega' := \bigcup \limits _{i=1} ^{m-1} \Omega_i$ and for all $i = 1,\dots,m-1$ let $\Omega_i' := \Omega_i \cap \Omega_m$. Then $(\Omega_i)_{i=1,\dots,m-1}$ is a reasonable cover of $\Omega'$ and $(\Omega_i')_{i =1,\dots,m-1}$ is a reasonable cover of $\Omega' \cap \Omega_m$. The Mayer-Vietoris-Sequence \ref{MVS} shows that for all $k$ the complex 
\begin{align} \label{sequence}
H^{k-1}(\Omega' \cap \Omega_m) \rightarrow H^{k}(\Omega)  \rightarrow H^{k} (\Omega') \oplus H^{k}(\Omega_m)
\end{align}
is exact. By induction hypothesis both $H^{k} (\Omega' \cap \Omega_m)$ and $H^{k} (\Omega')$ are finite dimensional and by definition so is $H^{k}(\Omega_m)$. Then by exactness of (\ref{sequence}), $H^{k}(\Omega)$ is finite dimensional. \\
If $(\Omega_i)_{i=1,\dots,m}$ is a good cover, then so are $(\Omega_i)_{i=1,\dots,m-1}$ and $(\Omega_i')_{i =1,\dots,m-1}$. So for $k \geq m$, by induction hypothesis $H^{k}(\Omega') = 0$ (since $k \geq m-1$) and $H^{k-1}(\Omega' \cap \Omega_m) = 0$ (since $k-1 \geq m-1$). Since then also $k \geq 2$ we further have $H^{k}(\Omega_m) = 0$ and (\ref{sequence}) becomes
\begin{align*}
0 \rightarrow H^{k}(\Omega) \rightarrow 0 \oplus 0,
\end{align*}
which shows $H^{k}(\Omega) = 0$. 
\end{proof}
\end{lem}

\subsection{Polyhedral Stars}

In this subsection $\CC$ will be a polyhedral complex in the real vector space $\R^r$. We will introduce a special class of open subsets of $|\CC|$ and show that these sets are polyhedrally star shaped. 

\begin{defn} 
Let $\sigma \in \CC$. We denote by $\mathring \sigma$ the \textit{relative interior} of $\sigma$, which is just $\sigma$ without its proper faces. We define the \textit{polyhedral star of $\sigma$} to be $\Omega_\sigma := \bigcup \limits_{\tau \in \CC, \sigma \prec \tau} \mathring \tau$. 
\end{defn}

\begin{lem}
For $\sigma \in \CC$ the polyhedral star $\Omega_\sigma$ of $\sigma$ is an open neighbourhood of $ \mathring \sigma$ in $|\CC|$. 
\begin{proof}
Since $\sigma \prec \sigma$, we have $ \mathring \sigma \subset \Omega_\sigma$. Let $z \in \Omega_\sigma$. Let $B$ be an open neighbourhood of $z$ in $\R^r$ that only intersects polyhedra in $\CC$ that contain $z$. Then we have $B \cap |\CC| \subset \bigcup \limits _{\tau : z \in \tau}  \mathring \tau$ and since $z \in \Omega_\sigma$ there exists some $\nu \in \CC$ such that $z \in  \mathring \nu $ and $\sigma \prec \nu$. Now if $z \in \tau$, then $z \in \nu \cap \tau$, which is a face of both. But since $z \in  \mathring \nu$ this can not be a proper face of $\nu$, hence $\nu \cap \tau = \nu$. Thus we have $\nu \prec \tau$ and by transitivity $\sigma \prec \tau$. We have shown $\{ \tau \in \CC | z \in \tau \} \subset \{ \tau \in \CC | \sigma \prec \tau \}$. This shows in turn that $ \bigcup \limits_{\tau : z \in \tau} \mathring \tau \subset \bigcup \limits _{\tau: \sigma \prec \tau}  \mathring \tau$ and thus $B \cap |\CC| \subset \bigcup \limits _{\tau: \sigma \prec \tau} \mathring \tau = \Omega_\sigma$. Hence for every point $z \in \Omega_\sigma$, the set $\Omega_\sigma$ contains an open neighbourhood of $z$ in $|\CC|$, which shows that $\Omega_\sigma \subset |\CC|$ is an open set.
\end{proof}
\end{lem}

\begin{lem}
Let $\tau_1,\dots,\tau_n \in \CC$. Then the set of polyhedra in $\CC$ which contain all $\tau_i$ is either empty or has a unique minimal (i.e. smallest) element $\sigma_{\tau_1\dots \tau_n}$. Further we have $\bigcap \limits_{i=1} ^{n} \Omega_{\tau_i} = \Omega_{\sigma_{\tau_1\dots \tau_n}}$.
\begin{proof}
The first assertion is clear since the set of polyhedra which contain all $\tau_i$ is closed under intersection. 
The second part is straight from the definition, since 
\begin{align*}
\bigcap \limits_{i=1} ^{n} \Omega_{\tau_i} = \bigcup \limits _{\nu: \tau_i \prec \nu \forall i} \mathring \nu = \bigcup \limits_{\nu: \sigma_{\tau_1\dots,\tau_n} \prec \nu} \mathring \nu = \Omega_{\sigma_{\tau_1,\dots,\tau_n}}.
\end{align*} 
\end{proof}
\end{lem}

\begin{lem} \label{lemstarshaped}
Let $\sigma \in \CC$ and $z \in  \mathring \sigma$. Then $\Omega_\sigma$ is polyhedrally star shaped with respect to $z$. 
\begin{proof}
Let $\DD$ be the polyhedral complex whose maximal polyhedra are the maximal ones in $\CC$ that contain $\sigma$. Let $\tau \in \DD$ be maximal and $y \in \tau \cap \Omega_\sigma$. Then there exists $\nu$ such that $y \in \mathring \nu$ and $\sigma \prec \nu \prec \tau$. Then $[y , z ) \subset \mathring \nu$ and hence $[y , z] \subset \mathring \nu \cup \mathring \sigma \subset \Omega_\sigma \cap \tau$. This just means that $\tau \cap \Omega_\sigma$ is star shaped, hence $\Omega_\sigma$ is polyhedrally star shaped.
\end{proof}
\end{lem}

\subsection{A finite dimensionality result for the operator $d'$}

In this subsection we will use the results of the previous two subsections together with the Poincar\'e lemma to shows that the cohomology with respect to $d'$ of superforms on polyhedral complexes is finite dimensional. Note that by symmetry for all statements for $d'$ the corresponding statements are true for $d''$. Again $\CC$ will be a polyhedral complex in $\R^r$. 

\begin{defn}
An open subset $\Omega \subset |\CC|$ is called \textit{polyhedrally connected} if there exists a polyhedral complex $\DD$ such that $\Omega$  is an open subset of $|\DD|$ and such that for each maximal polyhedron $\sigma$ in $\DD$ the set $\sigma \cap \Omega$ is connected.
\end{defn}

\begin{lem}  \label{Lemfinitedimensional}
Let $\Omega \subset |\CC|$ a polyhedrally connected open subset. Then $H^{0,q}_{d'}(\Omega)$ is a finite dimensional real vector space for all $q$. 
\begin{proof}
Choose a complex $\DD$ such that $\Omega \subset |\DD|$ is an open set and for all maximal $\sigma \in \DD$, the set $\sigma \cap \Omega$ is connected. By definition a superform $\alpha \in A^{p,q}(\Omega)$ on a polyhedral complex is closed under $d'$ if and only if all its restrictions $\alpha|_{\sigma \cap \Omega}$ to maximal polyhedra of $\DD$ are closed under $d'$. Hence the injection $A^{0,q}(\Omega) \inj \bigoplus \limits_{\text{maximal } \sigma \in \DD} A^{0,q}(\sigma \cap \Omega)$ restricts to an injection $H^{0,q}_{d'}(\Omega) \inj \bigoplus \limits_{\text{maximal } \sigma \in \DD} H^{0,q}_{d'}(\sigma \cap \Omega)$. It is easy to see that since $\sigma \cap \Omega$ is connected we have $H^{0,q}_{d'}(\sigma \cap \Omega) = \Lambda^q \Linear_\sigma^*$, where $\Linear_\sigma^*$ denotes the dual of the linear space associated to $\sigma$. Hence the sum is finite dimensional and thus $H^{0,q}_{d'}(\Omega)$ is. 
\end{proof}
\end{lem}

\begin{satz} \label{Thmd'finite}
$H^{p,q}_{d'}(|\CC|)$ is finite dimensional for all $p,q \in \N_0$. 
\begin{proof}
Let $\tau_1,\dots,\tau_k$ be the minimal polyhedra of $\CC$. We claim that the family $(\Omega_{\tau_i})_{i= 1,\dots,k}$ is a good cover of $|\CC|$. Let therefore $z \in |\CC|$. Then $z$ is in the relative interior of some polyhedron $\sigma$ and there is $\tau_i$ such that $\tau_i \prec \sigma$. This means however that $z \in \Omega_{\tau_i}$. Hence we have a cover and Lemma \ref{lemstarshaped} together with the Poincar\'e lemma (Theorem \ref{PLemmaComplex}) and Lemma \ref{Lemfinitedimensional} (using that polyhedrally star shaped sets are polyhedrally connected) precisely shows that this is a good cover. Now Lemma \ref{lemgoodcover} shows our result. 
\end{proof}
\end{satz}

\subsection{A finite dimensionality result for the operator $d$}

We will use the results of the previous three subsections and the Poincar\'e lemma to prove finite dimensionality for the cohomology with respect to $d$ of superforms on polyhedral complexes. Again $\CC$ will be a polyhedral complex of dimension $n$ in $\R^r$.

\begin{kor} \label{dPLemmaComplex}
Let $\CC$ be a polyhedral complex and $\Omega \subset |\CC|$ a polyhedrally star shaped open subset. Let $\alpha \in A^k(\Omega)$ be a $d$-closed form. Then there exists $\beta \in A^{k-1}(\Omega)$ such that $ \alpha - d\beta \in A^{0,k}(\Omega)$ and such that $\alpha - d\beta$ is $d', d''$ and $d$-closed. If $k > \dim \CC$ then $\alpha$ is $d$-exact.
\begin{proof}
Write $\alpha = \alpha_0 + \alpha_1 + \dots + \alpha_{k}$ with $\alpha_i \in A^{k-i, i}(\Omega)$. Then the decomposition of $d \alpha \in A^{k+1}(\Omega) = \bigoplus \limits _{p+q = k+1} A^{p,q}(\Omega)$ is given by
\begin{align*}
d\alpha = d' \alpha_0 + (d''\alpha_0 + d' \alpha_1) + \dots + (d''\alpha_{k-1} + d'\alpha_{k}) + d''\alpha_k.
\end{align*}
Since those terms have different bidegrees each of them is zero. Therefore the statement is trivially true if $k=0$ and we may from now on assume $k >0$. 

We construct inductively for $i = 0,\dots,k-1$ forms $\beta_i \in A^{k-i-1,i}(\Omega)$ such that $\beta_{-1} = 0$ and $d' \beta_i = \alpha_i - d'' \beta_{i-1}$. Note therefore that $\alpha_i - d'' \beta_{i-1}$ is $d'$-closed for $i = 0,\dots,k$, since this is immediate for $i = 0$ and for $i = 1,\dots,k$ we have
\begin{align*}
d'(\alpha_i - d'' \beta_{i-1}) &= d'\alpha_i - d' d'' \beta_{i-1} \\
&= d' \alpha_i + d'' d' \beta_{i-1} \\
&= d'\alpha_i + d'' \alpha_{i-1} - d'' d'' \beta_{i-2} \\
&= d'\alpha_i + d'' \alpha_{i-1} = 0.
\end{align*}
Hence given $\beta_{i-1}$, Theorem \ref{PLemmaComplex} gives us $\beta_i \in A^{k-i-1, i}(\Omega)$ such that $d' \beta_i = \alpha_i - d'' \beta_{i-1}$ for $i = 0,\dots,k-1$. We define $\beta := \sum \limits _{i=0} ^{k-1} \beta_i \in A^{k-1}(\Omega)$. Then we have
\begin{align*}
\alpha - d\beta = \sum \limits_{i = 0}^{k-1} (\alpha_i - d'' \beta_{i-1} - d' \beta_i) + \alpha_k - d'' \beta_{k-1} =  \alpha_k - d''\beta_{k-1} \in A^{0,k}(\Omega).
\end{align*}
As shown above we have that $\alpha_k - d'' \beta_{k-1}$ is $d'$ closed, thus $\alpha - d \beta$ is. Since it is also $d$-closed, it is $d''$-closed. If $k > \dim \CC$, then $A^{0,k}(\Omega) = 0$ and hence $\alpha = d \beta$.
\end{proof}
\end{kor} 

\begin{lem} \label{lemd'd''}
Let $\Omega \subset |\CC|$ be an open subset. Let $\alpha \in A^{0,k}(\Omega)$ such that $d' \alpha = 0$. Then $d'' \alpha = 0$ and $d \alpha = 0$. 
\begin{proof}
It is sufficient to check this after a restriction to a polyhedron. Let $\sigma \in \CC$ and let $v_1,\dots,v_r$ be a basis of $\Linear_\sigma$. Then $\alpha|_\sigma = \sum \limits _{|J| = k} \alpha_I d''v_J$ and $d' \alpha|_\sigma = 0$ if and only if $\frac {\del \alpha_J} {\del v_i} = 0$ for all $i, J$. But then also $d'' \alpha|_\sigma = 0$. \\
Since $d'\alpha = 0$ and $d''\alpha = 0$ we have $d\alpha = 0$. 
\end{proof}
\end{lem}

\begin{kor} \label{Kordcohomology}
Let $\Omega \subset |\CC|$ be a polyhedrally star shaped open subset. Then there is a surjective map $H^{0,k}_{d'}(\Omega) \surj H_d^k(\Omega)$. In particular $H_d^{k}(\Omega)$ is finite dimensional for all $k \in \N_0$.
\begin{proof}
By Lemma \ref{lemd'd''} the inclusion $A^{0,k}(\Omega) \inj A^{k}(\Omega)$ induces $H^{0,k}_{d'}(\Omega) \rightarrow H^k_d(\Omega)$ (note that $H_{d'}^{0,k} = \ker (d'_{0,k})$). Now Corollary \ref{dPLemmaComplex} shows the surjectivity and \ref{Lemfinitedimensional} shows that $H^{0,k}_{d'}(\Omega)$ is finite dimensional, hence $H^{k}_d(\Omega)$ is.
\end{proof}
\end{kor}

\begin{satz} \label{Thmdfinite}
$H^k_d(|\CC|)$ is finite dimensional for all $k \in \N_0$. 
\begin{proof}
Let $\tau_1,\dots,\tau_s$ be the minimal polyhedra in $\CC$. Again, as in the proof of \ref{Thmd'finite}, $(\Omega_{\tau_i})_{i=1,\dots,s}$ is a cover of $|\CC|$. By \ref{Kordcohomology} and \ref{lemstarshaped} this is a reasonable cover. Hence \ref{lemgoodcover} shows our result.
\end{proof}
\end{satz}

\section{Real-valued differential forms on Berkovich spaces}  \label{secPLemmaComplex}

In this section $K$ is a field which is algebraically closed and complete with respect to an absolute value. We work with a variety $X$ over $K$, by which we mean a reduced irreducible separated $K$-scheme of finite type. We let $n:= \dim(X)$ and denote by $\Xan$ the Berkovich analytification of $X$. \\
The space of real-valued $(p,q)$-forms on Berkovich analytic spaces was introduced by Chambert-Loir and Ducros in \cite{CLD} using analytic moment maps. In \cite{Gubler} Gubler developed an approach based on algebraic moment maps in the case where the analytic space is the analytification of an algebraic variety. In that case both approaches lead to the same sheaves of forms. We will follow Gubler's approach.

\subsection{A $d'$-Poincar\'e lemma for forms on Berkovich spaces}
\begin{defn}
\begin{enumerate} [itemindent =*, leftmargin=0mm]
\item
An open affine subset $U$ of $X$ is called \textit{very affine}, if it has a closed embedding to a torus $\G^s_m$, or equivalently if $\0_X(U)$ is generated by its units as a $K$-algebra. In this case $U$ has a canonical embedding (up to translation) $\varphi_U$ into a torus $\T_U$.  This embedding is constructed by choosing representatives of a basis $\varphi_1,\dots,\varphi_r$ of the free abelian group $M_U := \0_X(U)^\times / K^\times$, which yields a map $\varphi_U: U \rightarrow \T_U := \Spec K[M_U]$. The map $\varphi_U$ is called the \textit{canonical moment map} of $U$. We define the \textit{tropical variety} $\Trop(U)$ associated to $U$ to be the image of $\trop_U := \trop \circ (\varphi_U)^{\an}: U^{\an} \rightarrow N_{U,\R}$, where $N_{U} := {M_U}^*$, $N_{U,\R} := N_U \otimes_\Z \R$, $\trop: \T_U^{\an} \rightarrow N_{U,\R}$ is the tropicalization map of the torus and $(\varphi_U)^{\an}: U^{\an} \rightarrow \T^{\an}$ is the analytification of $\varphi_U$. It turns out that $\Trop(U)$ is the support of a polyhedral complex of pure dimension $n = \dim(X)$ in the $r$ dimensional real vector space $N_{U,\R}$. (cf. \cite[Theorem 3.3]{Gubler2})
\item
A \textit{tropical chart} is a pair $(V, \varphi_U)$, where $V \subset \Xan$ is an open subset in the analytic topology and $\varphi_U$ is the canonical moment map of a very affine Zariski open subset $U \subset X$, such that $V \subset U^{\an}$ and $V = \trop_U^{-1}(\Omega)$ for an open subset $\Omega$ of $\Trop(U)$.
\item
If $(V, \varphi_U)$ and $(V', \varphi_{U'})$ are tropical charts such that $V' \subset V$ and $U' \subset U$, then $(V', \varphi_{U'})$ is called a \textit{subchart} of $(V, \varphi_{U})$.
\end{enumerate}
\end{defn}

\begin{bem}
\begin{enumerate} [itemindent =*, leftmargin=0mm]
\item
The tropical variety $\Trop(U)$ actually has more structure. It is a \textit{rational} polyhedral complex (with respect to the lattice $N_U$), which is equipped with positive integer \textit{weights} on its top dimensional faces and satisfies the \textit{balancing condition} (cf. \cite{Gubler}). While these properties are used in the theory of differential forms, they are not needed for the aspects we consider.  
\item
For tropical charts $(V, \varphi_U)$ and $(V', \varphi_{U'})$ the pair $(V \cap V', \varphi_{U \cap U'})$ is a subchart of both. 
\item
Tropical charts form a basis of the topology of $\Xan$.
\item
If $(V', \varphi_{U'})$ is a subchart of $(V, \varphi_U)$ then there is a canonical surjective affine map 
\begin{align*}
\psi_{U', U}: N_{U',\R}\rightarrow N_{U, \R},
\end{align*}
with integral linear part, such that the diagram
\begin{align*}
\begin{xy}
\xymatrix
{
U' \ar[rr] \ar[d]_{\trop_{U'}} && U \ar[d]^{\trop_U} \\
\Trop(U') \ar[rr]_{\psi_{U',U}} && \Trop U
}
\end{xy}
\end{align*}
commutes.
\end{enumerate}
\end{bem}

\begin{defn}
Let $V$ be an open subset of $\Xan$. A $(p,q)$-differential form $\alpha$ on $V$ is given by a family $(V_i, \varphi_{U_i}, \alpha_i)_{i \in I}$ such that
\begin{enumerate} [itemindent =*, leftmargin=0mm]
\item
For all $i \in I$ the pair $(V_i, \varphi_{U_i})$ is a tropical chart and $\bigcup \limits _{i \in I} V_i = V$. 
\item
For all $i \in I$ we have $\alpha_i \in A^{p,q}(\trop_{U_i}(V_i))$.
\item
The $\alpha_i$ agree on intersections in the sense that for all $i,j \in I$, we have
\begin{align*}
\psi_{U_i \cap U_j, U_i}^*(\alpha_i) = \psi_{U_i \cap U_j, U_j}^*(\alpha_j) \in A^{p,q}(\trop_{U_i \cap U_j}(V_i \cap V_j)).
\end{align*}
\end{enumerate}
Another such family $(V'_j, \varphi_{U'_j}, \beta_j)_{j \in J}$ defines the same form if there is a common refinement of the covers of $V$ by tropical charts such that the affine pullbacks to the refined cover agree.\\
We write $A^{p,q}(V)$ for the space of differential forms of bidegree $(p,q)$ on $V$ and $A^{p,q}$ for the sheaf of differential forms of bidegree $(p,q)$ on $\Xan$. We also write $A^k := \bigoplus \limits _{p+q = k} A^{p,q}$ for the sheaf of differential forms of degree $k$. \\
Since affine pullbacks are compatible with $d'$, we can define $d'\alpha$ to be given by $(V_i, \varphi_{U_i}, d'\alpha_i)_{i \in I}$. This defines a well defined operator $d': A^{p,q}(V) \rightarrow A^{p+1,q}(V)$. The same works for $d''$ and $d$ and we get differential operators $d'': A^{p,q}(V) \rightarrow A^{p,q+1}(V)$ and $d: A^k(V) \rightarrow A^{k+1}(V)$. 
\end{defn}

\begin{bem}
It is obvious that $d'$ is a differential. Hence for each $q \in \{0,\dots,n\}$ we get a complex
\begin{align*}
0 \rightarrow A^{0,q} \overset {d'} {\rightarrow} A^{1,q} \overset {d'} {\rightarrow} \dots \overset {d'} {\rightarrow} A^{n,q} \rightarrow 0
\end{align*}
 of sheaves on $\Xan$. Theorem \ref{PLemmaXan} will show that this complex is always exact in positive degrees. We also get a complex 
\begin{align*}
0 \rightarrow A^{0} \overset {d}{\rightarrow} A^1 \overset{d}{\rightarrow} \dots \overset {d}{\rightarrow} A^{2n} \rightarrow 0
\end{align*}
of sheaves on $X$, but we can not hope for this complex to be exact for the same symmetry reason as given in \ref{remnotexact}.
\end{bem}

\begin{satz} [$d'$-Poincar\'e lemma on $\Xan$] \label{PLemmaXan}
Let $X$ be a variety and $V \subset \Xan$ an open subset. Let $x \in V$ and $\alpha \in A^{p,q}(V)$ with $p > 0$ and $d'\alpha = 0$. Then there exists some open $W \subset V$ with $x \in W$ and some $\beta \in A^{p-1,q}(W)$ such that $d'\beta = \alpha|_W$.
\begin{proof}
Let $\alpha$ be given by a family $(V_i, \varphi_{U_i}, \alpha_i)_{i \in I}$ where $(V_i, \varphi_{U_i})$ are tropical charts, $\alpha_i \in A^{p,q}(\Omega_i)$ and $\Omega_i := \trop_{U_i}(V_i)$ is an open subset of $\Trop(U_i)$. Choose $i$ such that $x \in V_i$ and let $z := \trop_{U_i}(x)$. By Remark \ref{bemlocallystarshaped} we may choose a polyhedrally star shaped neighbourhood $\Omega'$ of $z$ in $\Omega_i$. We define $W := \trop_{U_i}^{-1} (\Omega')$. Then $\alpha|_W$ is given by the single chart $(W, \varphi_{U_i}, \alpha_i|_{\Omega'})$ and $d'\alpha|_W$ is given by $(W, \varphi_{U_i}, d'\alpha_i|_{\Omega'})$. Since $d'\alpha|_W = 0$ and it is given by a single chart, we know that $d'\alpha_i|_{\Omega'} = 0$ \cite[Proposition 5.6]{Gubler}. Now Theorem $\ref{PLemmaComplex}$ applies and gives us $\beta' \in A^{p-1,q}(\Omega')$ such that $d'\beta' = \alpha_i|_{\Omega'}$. The form $\beta \in A^{p-1,q}(W)$ given by $(W, \varphi_{U_i}, \beta')$ now has the desired property.
\end{proof}
\end{satz}

\begin{kor} \label{korPLemmaXan}
The complex
\begin{align} \label{acycresspace}
0 \rightarrow \underline{\R} \rightarrow A^{0,0} \overset{d'} {\rightarrow} A^{1,0} \overset {d'} {\rightarrow} \dots \overset{d'} {\rightarrow}  A^{n,0} \rightarrow 0
\end{align}
of sheaves on $\Xan$ is exact. The cohomology of its complex of global sections
\begin{align} \label{globsecspace}
0 \rightarrow A^{0,0}(\Xan) \overset{d'} {\rightarrow} A^{1,0}(\Xan) \overset {d'} {\rightarrow} \dots \overset{d'} {\rightarrow}  A^{n,0}(\Xan) \rightarrow 0
\end{align}
is isomorphic to the sheaf cohomology $H^*(\Xan, \underline{\R})$ of the constant sheaf $\underline{\R}$ and to the singular cohomology $H_{\sing}^*(\Xan, \R)$.
\begin{proof}
Let $V \subset \Xan$ an open subset and $f \in A^{0,0}(V)$ given by the family $(V_i, \varphi_{U_i}, f_i)_{i \in I}$. Then $f$ can be viewed as a function on $V$ via the definition $f(x) := f_i \circ \trop_{U_i}(x)$ for $x \in U^{\an}_i$. This function is continuous and $d'f = 0$ if and only if $f$ is locally constant. Together with Theorem $\ref{PLemmaXan}$ this shows that the complex (\ref{acycresspace}) is exact. Since the sheaves $A^{p,q}$ admit partitions of unity \cite[Proposition 5.10]{Gubler}, they are fine, hence acyclic \cite[Chapter II, Proposition 3.5 and Theorem 3.11]{Wells}. This means that the complex (\ref{globsecspace}) calculates the sheaf cohomology of $\underline{\R}$. Since $\Xan$ is paracompact, Hausdorff and locally compact, this is the singular cohomology of the underling topological space of $\Xan$ \cite[Chapter III, Theorem 1.1]{Bredon}.
\end{proof}
\end{kor}
 
\begin{kor}
Let $X$ be a variety and $V \subset \Xan$ an open subset. Let $\alpha \in A^{k}(V)$ such that $d\alpha = 0$. Then for $x \in V$ there exists an open neighbourhood $W$ of $x$ in $V$ and a form $\beta \in A^{k-1}(W)$ such that $\alpha|_W - d\beta \in A^{0,k}(W)$ and such that $\alpha|_W - d\beta$ is closed under $d$, $d'$ and $d''$. If $k > \dim(X)$ then $\alpha|_W$ is $d$-exact.
\begin{proof}
The proof works the same as the proof of Theorem \ref{PLemmaXan}, using \ref{dPLemmaComplex} instead of \ref{PLemmaComplex}.
\end{proof}
\end{kor}

\begin{satz} \label{PLemmaCLD}
Let $X$ be a Berkovich analytic space of dimension $n$. Let $A^{p,q}$ be the sheaf of differential $(p,q)$-forms on X as introduced by Chambert-Loir and Ducros in \cite{CLD}. Then for all $q \in \{0,\dots,n\}$ the complex 
\begin{align*}
0 \rightarrow  A^{0,q} \overset{d'} {\rightarrow} A^{1,q} \overset {d'} {\rightarrow} \dots \overset{d'} {\rightarrow}  A^{n,q} \rightarrow 0
\end{align*}
of sheaves on $X$ is exact in positive degrees. Further the complex
\begin{align*}
0 \rightarrow \underline{\R} \rightarrow A^{0,0} \overset{d'} {\rightarrow} A^{1,0} \overset {d'} {\rightarrow} \dots \overset{d'} {\rightarrow}  A^{n,0} \rightarrow 0
\end{align*}
of sheaves on $X$ is exact. \\
If $X$ is a good analytic space which is Hausdorff and paracompact, then the cohomology of the complex 
\begin{align*}
0 \rightarrow A^{0,0}(X) \overset{d'} {\rightarrow} A^{1,0}(X) \overset {d'} {\rightarrow} \dots \overset{d'} {\rightarrow}  A^{n,0}(X) \rightarrow 0
\end{align*}
is equal to the sheaf cohomology $H^*(X, \underline{\R})$ of the constant sheaf $\underline{\R}$, which is isomorphic to the singular cohomology $H_{\sing}^*(X, \R)$. 
\begin{proof}
Using \cite[Lemme 3.2.2]{CLD} the same arguments as used in the proof of Theorem \ref{PLemmaXan} work, since forms in the sense of \cite{CLD} are also locally given by forms on polyhedral complexes. If $X$ is good, Hausdorff and paracompact, then \cite[Proposition 3.3.6]{CLD} shows that there are partitions of unity and the arguments in the proof of Corollary \ref{korPLemmaXan} work. The details are left to the reader.
\end{proof}
\end{satz}

As observed in Remark \ref{bemsymmetry} the corresponding statements to Theorem \ref{PLemmaXan} and \ref{korPLemmaXan} are true for $d''$.

\subsection{Finite dimensionality of the de Rham cohomology of differential forms on Berkovich spaces} \label{ssecfindim}
We will now again work with the case where our analytic space is the analytification of an algebraic variety $X$. We first note that Theorem \ref{MVS} and Lemma \ref{lemgoodcover} apply also for forms on $X$ and open covers of $X$, since the proofs work exactly the same. However we only have Theorem \ref{PLemmaXan} available, which does not tell us anything about acyclic domains and hence we can not use the same strategy as in Section \ref{secfindim} to proof finite dimensionality results. However Theorem \ref{PLemmaCLD} shows that the cohomology of the complex $(A^{\bullet,0}(\Xan),d')$ is equal to singular cohomology, hence only depends on the homotopy type of the underlying topological space.  This gives us the possibility to prove finite dimensionality results for this complex.

\begin{satz} \label{ThmSkeleton}
Let $X$ be a variety. Then $H^{p,0}_{d'}(\Xan)$ is finite dimensional for all $p$.
\begin{proof}
We start with the quasi-projective case. By \cite[Theorem 13.2.1]{Skeletons} there exists a strong deformation retraction of $\Xan$ to a finite simplicial complex $S$. Finite simplicial complexes have finite dimensional singular cohomology as is certainly well known from algebraic topology. Since we also have $H^{p,0}_{d'}(\Xan) = H^p_{\sing}(\Xan, \R) = H^p_{\sing}(S, \R)$ the result follows for quasi-projective varieties. In the general case, let $X = \bigcup_{i=1} ^k U_i$ be a cover of $X$ by affine schemes. Since $X$ is separated, all intersections of the $U_i$ are affine and since affine schemes are quasi-projective, this is a reasonable cover. This shows our result.  
\end{proof}
\end{satz}

\begin{bem}
There are more result, based on the existence of (formal) models, that state that a given Berkovich space admits a deformation retraction to a so called skeleton, which is a finite polyhedral complex. Since this is the only property which we need in the proof of the previous theorem, we get corresponding finite dimensionality results in all these cases as well. 
\end{bem}

\bibliographystyle{alpha}

\small{Philipp Jell, Fakult\"at f\"ur Mathematik, Universit\"at Regensburg \\ Universit\"atsstra\ss e 31, 93053 Regensburg, Germany \\ Philipp.Jell@mathematik.uni-regensburg.de}

\end{document}